\magnification=\magstep1
\parskip=0mm
\hsize=6.25 true in
\frenchspacing
\bigskipamount=24pt plus 6 pt minus 6 pt
\medskipamount=18pt plus 3 pt minus 3 pt
\smallskipamount=12pt plus 2 pt minus 2 pt
\thinmuskip=5mu
\medmuskip=6mu plus 3mu minus 2mu
\thickmuskip=6mu plus 3mu minus 2mu
\def\7{{\hskip-4pt}}
\overfullrule=0pt

\font\xivrm=cmr17


%

\def\R{{\rm I\mkern-4mu R}}
\def\N{{\rm I\mkern-4mu N}}

\def\Z{{\rm Z\mkern-6mu Z}}

\def\Any{\,\cdot\,}

\def\modul#1{\mathopen\vert #1 \mathclose\vert}

\def\mod{\mathop{\rm mod }}

\def\Dd#1#2{{d^{#2} \over {d #1}^{#2}}}

\def\giantrightarrow{\mathrel{\hbox to 1.5cm{\rightarrowfill}}}
\def\ghookrightarrow{\mathrel{\hbox to 1.5cm{$\lhook\mkern-9mu$\rightarrowfill}}}

\def\teta{\vartheta}
\def\ro{\varrho}
\def\pfi{\varphi}
\def\Thm#1{\par\medskip\noindent{\bf #1}\qquad\bgroup\it}
\def\endThm{\egroup\par\medskip\noindent}

\centerline{\xivrm Exceptional Coupling Constants for the Coulomb-Dirac}
\centerline{\xivrm Operator with Anomalous Magnetic Moment}

\bigskip
\centerline{K.M. Schmidt}
\smallskip
\centerline{\it School of Mathematics, Cardiff University,\/}
\centerline{\it	P.O. Box 926, Cardiff CF24 4YH, UK\/}

\bigskip
\centerline{\it Dedicated to E.B. Davies on the occasion of his 60th birthday\/}

\bigskip
{\narrower\noindent
{\bf Abstract.}\quad
It was recently shown that the point spectrum of the separated Coulomb-Dirac
operator $H_0(k)$ is the limit of the point spectrum of the Dirac operator
with anomalous magnetic moment $H_a(k)$ as the anomaly parameter tends to $0$;
this spectral stability holds for all Coulomb coupling constants $c$ for which
$H_0(k)$ has a distinguished self-adjoint extension if the angular momentum
quantum number $k$ is negative, but for positive $k$ there are certain
exceptional values for $c$.
Here we obtain an explicit formula for these exceptional values.
In particular, it implies spectral stability for the three-dimensional
Coulomb-Dirac operator if $\modul{c} < 1$, covering all physically relevant
cases.
}

\bigskip
\noindent
{\bf 1. Introduction}

\medskip
\noindent
By separation of variables in spherical polar coordinates, the Dirac
operator of relativistic quantum mechanics with a Coulomb potential,
$$
	{\cal H}_0 := -i \alpha \cdot \nabla + \beta + {c \over \modul{x}},
$$
where $c < 0$ and $\alpha_1, \alpha_2, \alpha_3$ and $\beta = \alpha_0$
are symmetric $4\times 4$ matrices satisfying the anticommutation
relations
$$
	\alpha_i \alpha_j + \alpha_j \alpha_i = \delta_{ij}
\qquad	(i, j \in \{0, 1, 2, 3\}),
$$
is unitarily equivalent to a direct sum of one-dimensional Dirac operators
on the half-line,
$$
	H_0(k) = -i \sigma_2 \Dd{r}{} + \sigma_3 + {k \over r} \sigma_1
		+ {c \over r}
\qquad	(r \in (0, \infty)),
$$
where $\sigma_1 = \pmatrix{0 & 1 \cr 1 & 0 \cr}$,
$\sigma_2 = \pmatrix{0 & -i \cr i & 0 \cr}$ and
$\sigma_3 = \pmatrix{1 & 0 \cr 0 & -1 \cr}$ are the Pauli matrices and
$k \in \Z \setminus \{0\}$ is the angular momentum quantum number [13,
Appendix to Section 1].
For an electron orbiting a nucleus of charge number $Z \in \N$, the
Coulomb coupling constant is $c = - Z \alpha$, with the Sommerfeld fine
structure constant $\alpha \approx 1/137$.

	$H_0(k)$ is essentially self-adjoint on its minimal domain if and
only if $c^2 \le k^2 - 1/4$ [13, Theorem 6.9].
For $c^2 \in (k^2 - 1/4, k^2)$, there still exists a distinguished self-adjoint
extension, characterised by the requirement that the wave-functions in its
domain have a finite potential (or kinetic) energy [13, Theorem 6.10].

	Pauli suggested a modification of the Dirac operator
which takes into account the anomalous magnetic moment of the electron
(for the historical background cf. [7]). With suitably normalised constants,
the operator with a potential $V$ for an electron of magnetic moment
$(1 + a) \mu_B$, where $\mu_B$ is the Bohr magneton, takes the form
$$
	{\cal H}_a = -i \alpha \cdot \nabla + \beta + V
		 - i a \alpha \cdot \nabla V.
$$
In the case of the Coulomb potential, the corresponding half-line operators
after variable separation will be
$$
	H_a(k) = -i \sigma_2 \Dd{r}{} + \sigma_3
		+ \left({k \over r} + {a \over r^2}\right) \sigma_1
		+ {c \over r}
\qquad	(r \in (0, \infty)).
$$
The mathematical investigation of the properties of this operator was
initiated by Behncke [2, 3, 4].
$H_a(k)$ is essentially self-adjoint on its minimal domain for all values
$c < 0$ of the coupling constant (cf. [1, 5]).
The essential spectrum of $H_a(k)$, as that of $H_0(k)$, is
$(-\infty, -1] \cup [1, \infty)$.

	Let us now focus on the case $c^2 < k^2$, so that a distinguished
self-adjoint realisation of $H_0(k)$ exists.
We always assume $a < 0$ in the following; the case of positive $a$ can be
reduced to this (with a sign change of $k$) by a suitable unitary
transformation.
Although only integer values
of $k$ are relevant for the three-dimensional operator, we admit general
non-zero real values of $k$.
As is well known, $H_0(k)$ has infinitely many eigenvalues in the spectral
gap $(-1, 1)$, which accumulate at $1$. One would
expect that the eigenvalues of $H_a(k)$ will be perturbations of those of
$H_0(k)$, such that each eigenvalue of $H_0(k)$ will be the limit of exactly
one eigenvalue branch of $H_a(k)$ in the limit $a \rightarrow 0$
({\it spectral stability\/}, cf. [8, Chapter VIII \S 1.4]).
This expectation is partly corroborated by the strong resolvent convergence
of $H_a(k)$ to $H_0(k)$, at least for $c^2 \le k^2 - 1/4$, which implies that
the spectrum cannot suddenly expand in the limit.
Nevertheless, due to the strong singularity of the $a/r^2$ term at the
origin, this limit problem is highly non-trivial and does not subject itself
to the standard perturbation techniques.
Behncke [4, Theorem 2] was able to prove spectral stability for
$c^2 < k^2 - (3/2)^2$ if $k < 0$, and for $c^2 < k - 5/2$ if $k > 0$.

	The surprising asymmetry with respect to the sign of $k$ is not an
artefact of Behncke's method, but inherent in the problem, as observed in
the recent treatment [7],
where spectral stability was proved, based on an asymptotic analysis of the
(non-linear) Pr\"ufer and Riccati equations equivalent to the Dirac system
$$
	(H_a(k) - \lambda) u = 0,
\eqno	(1)
$$
for all $c \in (k, 0)$ if $k < 0$, and for all $c \in (-k, 0) \setminus
\{c_0, c_1, \dots\}$ if $k > 0$.
Here $c_j \in (-k, 0)$, with $c_j > c_{j+1}$, are certain exceptional values
of the coupling constant at which a transition in the behaviour of the
eigenfunctions of $H_a(k)$ occurs: for $c \in (c_n, c_{n-1})$, the eigenfunctions
of $H_a(k)$ show $n$ additional rapid oscillations very close to the origin,
compared to the corresponding eigenfunctions of $H_0(k)$.
It remains a fairly subtle open question whether or not spectral stability
holds if $c$ is equal to one of the exceptional values.

	The present note is devoted to a study of the number and positions
of the exceptional values of $c$ for a given $k > 0$.

\medskip
\noindent
{\bf 2. Exceptional values}

\medskip
\noindent
The exceptional values for $c$ appear in a stability
analysis of the differential equation [7, Eq.(1)]
$$
	\ro \teta'(\ro) = c + (k - {1 \over \ro}) \sin 2 \teta
\qquad	(\ro > 0).
\eqno	(2)
$$
This equation arises from the Pr\"ufer transformation (cf. Appendix A) of the
Dirac system $(1)$ after omitting the lowest-order term at $0$,
$\sigma_3 - \lambda$, and rescaling $r = \modul{a}\ro$ to absorb the parameter
$a$.

	In the limit $\ro \rightarrow \infty$, the right-hand side of
$(2)$ --- which is $\pi$-periodic in $\teta$ --- has asymptotic zeros
$\teta_\pm(c, k)$ satisfying
$0 < \teta_-(c, k) < \pi/4 < \teta_+(c, k) < \pi/2$,
$$
	\sin 2 \teta_\pm(c, k) = - {c \over k},
\qquad	\tan \teta_\pm(c, k) = {k \pm \sqrt{k^2 - c^2} \over -c}.
$$
An asymptotic study of the direction field of $(2)$ for $\ro \rightarrow 0$
and $\ro \rightarrow \infty$ yields the following result, which shows that
this differential equation has (up to addition of a constant multiple of $\pi$)
a unique unstable solution at $0$ and at $\infty$ (see [7, Lemma 2.1, 2.2]).
\Thm{Lemma 1.}
Let $k > 0$ and $c\in (-k, 0)$. Then the following statements hold.

\item{a)}
There is a unique solution $\teta_0(\Any, c)$ of $(2)$ such that
$\lim\limits_{\ro \rightarrow 0} \teta_0(\ro, c) = \pi$.

\item{}
All other solutions either differ from $\teta_0$ by a constant integer multiple
of $\pi$, or else satisfy
$\lim\limits_{\ro \rightarrow 0} \teta(\ro, c) = {\pi \over 2} \mod \pi$.

\item{}
For fixed $\hat\ro > 0$, $\teta_0(\hat\ro, \Any)$ is continuous non-decreasing.

\item{b)}
There is a unique solution $\teta_\infty(\Any, c)$ of $(2)$ such that
$\lim\limits_{\ro \rightarrow \infty} \teta_\infty(\ro, c) = \teta_-(c, k)$.

\item{}
All other solutions either differ from $\teta_\infty$ by a constant integer
multiple of $\pi$, or else satisfy
$\lim\limits_{\ro \rightarrow \infty} \teta(\ro, c) = \teta_+(c, k) \mod \pi$.

\item{}
For fixed $\hat\ro > 0$, $\teta_\infty(\hat\ro, \Any)$ is continuous and
strictly decreasing.

\endThm
The solution $\teta_0(\Any, c)$ asymptotically corresponds to the Pr\"ufer angle
of the $L^2(0, \infty)$ solution of $(1)$.

	The {\it exceptional values\/} are defined as those values of $c$ for
which $\teta_0(\Any, c)$ and $\teta_\infty(\Any, c)$ match up $\mod \pi$,
so that the unstable solution of $(2)$ at $0$ is unstable at infinity as well.
More precisely, in view of the monotonicity properties of $\teta_0$
and $\teta_\infty$ with respect to $c$, we have
$$
	\lim_{\ro \rightarrow \infty} \teta_0(\ro, c_m)
	= \teta_-(c_m, k) - m \pi
\qquad	(m = 0, 1, 2, \dots).
$$

	By the transformation $\ro = e^t$, $\teta(\ro) = \pfi(\log\ro)$,
the equation $(2)$ is equivalent to
$$
	\pfi'(t) = c + (k - e^{-t}) \sin 2 \pfi(t)
\qquad	(t \in \R),
\eqno	(3)
$$
the differential equation for the Pr\"ufer angle
$\pfi$ (cf. Appendix A) of a $\R^2$-valued solution $u$ of the Dirac
system
$$
	(-i \sigma_2 \Dd{t}{} + (k - e^{-t}) \sigma_1 + c)\, u(t) = 0,
\eqno	(4)
$$
which can be rewritten using the definition of the Pauli matrices as
$$	\eqalign{
	u_1' &= (e^{-t} - k) u_1 - c u_2,
\cr	u_2' &= c u_1 + (k - e^{-t}) u_2.
\cr	}
$$

	In $(4)$, $c$ takes on the role of a spectral parameter, while the
coefficient of $\sigma_1$ can be interpreted as a constant mass term with
an exponentially decaying perturbation. In the sense of the analogue of
Kneser's Theorem for Dirac systems [9], this perturbation is subcritical,
and indeed the method developed in [9] can be used to show that, for each
$k > 0$, there are only finitely many exceptional values for $c$.
Moreover, an asymptotic analysis of $(4)$ along the lines of
[10; see also 12] reveals that the number of exceptional
values is asymptotic to $k$ in the limit $k \rightarrow \infty$ with fixed
ratio $c / k$.

	This was verified in a computational investigation of $(4)$
based on a piecewise-constant approximation of the exponential function,
following the approach of [11].
The numerical findings suggested that the number of exceptional values
in $(-k, 0)$ is always equal to the unique integer in $[k - 1, k)$; more
precisely, whenever $k$ reaches an integer value, an additional exceptional
value, initially equal to $-k$, appears and moves into the interval $(-k ,0)$
as $k$ increases.

	This regularity raised the suspicion that $(4)$ has an underlying
solvable structure, and indeed, on closer scrutiny it turns out that this
equation
can be analysed by means of a variant of the factorisation method,
resulting not only in a proof of the above conjecture, but even in an
explicit formula for the exceptional values, which makes any asymptotic and
computational analysis of (4) superfluous.
\Thm{Theorem 1.}
For given $k > 0$, the exceptional values for $H_a(k)$ are given by
$$
	c_{n-1}(k) = - \sqrt{2 k n - n^2}
\qquad	(n \in \N, n < k).
$$
Consequently there are exactly $N$ exceptional values in $(-k, 0)$ if
$k \in (N, N + 1]$, $N \in \N_0$.
\endThm
In view of Theorems 1.1 and 1.2 of [7], this implies in particular that the
eigenvalues of $H_a(k)$ converge to those of $H_0(k)$ as $a \rightarrow 0$
if $c^2 < 2 k - 1$, confirming Behncke's linear bound (although with
a different constant). However, his conjecture [4, p. 2558] that spectral
stability may hold for all $c^2 < k^2 - 5/2$ turns out to be very questionable.

	Since the three-dimensional Coulomb-Dirac Hamiltonian ${\cal H}_a$
is the direct sum of $H_a(k)$ for $k \in \Z \setminus \{0\}$, we can thus infer
spectral stability for this operator in all cases in which ${\cal H}_0$ has a
distinguished self-adjoint realisation.
\Thm{Corollary 1.}
Let ${\cal H}_0$ be the Dirac operator and ${\cal H}_a$ the Dirac operator
with anomalous magnetic moment with a Coulomb potential with coupling constant
$c \in (-1, 0)$.
Then the eigenvalues of ${\cal H}_0$ are the limits of the
eigenvalues of ${\cal H}_a$ as $a \rightarrow 0$.
\endThm
{\bf 3. Proof of Theorem 1}

\medskip
\noindent
The following proof of Theorem 1 will be based on squaring the Dirac-type
operator in $(4)$, which is equivalent to deriving second-order differential
equations for $u_1$ and $u_2$.
In a general situation, this is usually not a good idea except to obtain a
quick and rough heuristic result, since it involves derivatives of the
coefficients and consequently requires unnecessarily strong regularity.
In the present instance, however, we are dealing with a specific Dirac
system with analytic coefficients, and as it turns out, the resulting
second-order equation system decouples and can be solved by the factorisation
method in the relevant cases (the second-order equation (5) is closely
related to the radial Schr\"odinger equation with a Morse potential [6, Section
5.2]). Here the factorisation is incidentally
provided by the original Dirac system; see the discussion in Appendix B.
The key observation is contained in the following result.
\Thm{Theorem 2.}
Let $\kappa \ge 0$. Let $v_1(t, \kappa) := e^{-e^{-t} - \kappa t}$,
and define $v_j$ recursively for $j \in \N + 1$ by
$$
	v_{j + 1}(t, \kappa) := (\kappa + j - e^{-t}) v_j(t) - v_j'(t)
\qquad	(t \in \R).
$$
Then $v_j(\Any, \kappa)$ is a nontrivial solution of
$$
	v'' = (e^{-2t} - (2 \kappa + 2 j - 1) e^{-t} + \kappa^2) \,v,
\eqno	(5)
$$
and it has the asymptotic properties
$$
	\lim_{t \rightarrow \infty} v_j(t,\kappa) = 0,
\qquad	\lim_{t \rightarrow \infty} {v_j'(t,\kappa) \over v_j(t,\kappa)}
		= - \kappa,
\quad	\hbox{and}
\quad	\lim_{t \rightarrow -\infty} {v_j'(t, \kappa) \over v_j(t, \kappa)}
		e^t = 1
\qquad	(j \in \N).
$$
\/}

\medskip
\noindent
The {\it proof\/} can be done by induction with respect to $j$.
The fact that $v_j$ is a solution of $(5)$ is readily verified by
differentiating twice.
The asymptotic properties are obvious for
$$
	{v_1'(t, \kappa) \over v_1(t, \kappa)} = e^{-t} - \kappa.
$$

	Now assume that $j \in \N$ is such that the assertion is true.
If $v_{j+1}$ were trivial, this would imply that (up to multiplication by
a constant) $v_j(t, \kappa) = e^{(\kappa + j)t + e^{-t}} \rightarrow \infty$
$(t \rightarrow \infty)$, contradicting the first limit property of $v_j$.
The asymptotic properties of $v_{j+1}$ are easily checked using the
identity
$$
	{v_{j+1}'(t, \kappa) \over v_{j+1}(t, \kappa)}
	= {-e^{-2t} + 2(\kappa+j) e^{-t} - \kappa^2) + (\kappa + j - e^{-t})
		{v_j'\over v_j}(t, \kappa)
\over	\kappa + j - e^{-t} - {v_j'\over v_j}(t, \kappa) }.
$$
\Thm{Corollary 2.}
Let $\kappa \ge 0$, $j \in \N$ and $v_j$ be as in Theorem 1. Then
$$
	\lim_{t \rightarrow \infty} {v_j(t) \over v_{j+1}(t)}
		= {1 \over 2 \kappa + j},
\qquad	\lim_{t \rightarrow -\infty} {v_j(t) \over v_{j+1}(t)} = 0.
$$
\endThm
This is a direct consequence of Theorem 1 in view of
$$
	{v_j(t, \kappa) \over v_{j+1}(t, \kappa)}
		= {1 \over \kappa + j - e^{-t} - {v_j'\over v_j}(t, \kappa) }.
$$
Furthermore, the following conclusion can be verified by a straightforward
calculation.
\Thm{Corollary 3.}
Let $v_j$ be defined as in Theorem 1, and let $j \in \N$, $ k \ge j$.
Then
$$
	u(t) = \pmatrix{
		{1 \over \sqrt{2 k j - j^2}} v_{j+1}(t, k - j)	\cr
		v_j(t, k - j)					\cr}
\qquad	(t \in \R)
$$
is a nontrivial solution of the Dirac system $(4)$ with
$c = -\sqrt{2 k j - j^2}$.
The Pr\"ufer angle $\pfi$ of $u$ satisfies
$$
	\lim_{t \rightarrow \infty} \tan \pfi(t) = {j \over \sqrt{2kj - j^2}},
\qquad	\lim_{t \rightarrow -\infty} \tan \pfi(t) = 0.
$$
\endThm
Since, with $c$ as in Corollary 3,
$$
	\tan \teta_-(c) = {k - \sqrt{k^2 - c^2} \over -c}
	= {k - \sqrt{k^2 - 2 k j + j^2} \over \sqrt{2 k j - j^2}}
	= {j \over \sqrt{2 k j - j^2}},
$$
where we have used $j \le k$, the Pr\"ufer angle $\pfi$ in Corollary 3
corresponds via $\teta(\ro) = \pfi(\log\ro)$ to a solution of $(3)$ which
coincides (up to addition of integer multiples of $\pi$) with both
$\vartheta_0(\Any, c)$ and $\vartheta_\infty(\Any, c)$; hence such values
of $c$ are exceptional values.

	We now conclude the proof of Theorem 1 by showing, based on the
continuity and monotonicity properties of $\vartheta_0$ and $\vartheta_\infty$,
that all exceptional values are obtained in this way. Corollary 3 specifies
the asymptotic behaviour of the relevant angle functions $\mod \pi$;
a study of the zeros of $v_j$ reveals their global behaviour as shown in
Lemma 2 below. Theorem 1 then follows in view of
$$
	\teta(\ro, -\sqrt{2 k j - j^2}) = \pfi_j(\log \ro, k - j)
\qquad	(\ro > 0, k > 0, j \in \N, j < k),
$$
where (with $\arctan 0 = 0$)
$$
	\pfi_j(\Any, \kappa)
		:= \arctan ({v_j(\Any, \kappa) \over v_{j+1}(\Any, \kappa)}
		\sqrt{2 \kappa j + j^2}) + \pi.
$$
\Thm{Lemma 2.}
Let $\kappa \ge 0$. Then $v_j(\Any, \kappa)$ has exactly $j - 1$ zeros.
The associated angle function $\pfi_j$ satisfies
$$
	\lim_{t \rightarrow -\infty} \pfi_j(t, \kappa) = \pi
$$
and
$$
	\lim_{t \rightarrow \infty} \pfi_j(t, \kappa)
		= \teta_-(-\sqrt{2 \kappa j + j^2}, \kappa + j) - (j-1)\pi
\qquad	(j \in \N).
$$
\endThm
{\it Proof.\/}
A look at the direction field of $(3)$ with $k = \kappa + j$ shows that
$\pfi_j(\Any, \kappa)$,
being the unstable solution at both $\pm \infty$, is monotone decreasing.
$v_1(\Any, \kappa)$ has no zeros, so $\pfi_1(\Any, \kappa)$ has no zeros
$\mod \pi$, and the assertion follows for $j = 1$.

	To conclude the proof by induction, assume now that $j \in \N$ is
such that the assertion is true.
Then $\pfi_j(\Any, \kappa) - \pi/2$ has exactly $j$ zeros mod $\pi$,
so $v_{j+1}(\Any, \kappa)$ has
exactly $j$ zeros.
Since the $v_{j+1}(\Any, \kappa)$ is a non-trivial solution
of a second-order equation,
$v_{j+1}'(\Any,\kappa)$ does not vanish at zeros of $v_{j+1}(\Any, \kappa)$.
Consequently,
$v_{j+1}(\Any, \kappa)$ and $v_{j+2}(\Any, \kappa)$
have no common zeros.
Hence $\pfi_{j+1}(\Any, \kappa)$ has exactly $j$ zeros $\mod \pi$,
and therefore must
converge to $\vartheta_-(- \sqrt{2 \kappa j + j^2, \kappa + j}) - (j-1)\pi$
as $t \rightarrow \infty$.

\medskip
\noindent
{\bf Appendix}

\medskip
\noindent
{\it A. The Pr\"ufer Transformation\/}

\noindent
Consider a general Dirac system
$$
	(-i \sigma_2 \Dd{x}{} + m(x) \sigma_3 + l(x) \sigma_1 + q(x)) u (x)
	= 0
\qquad	(x \in I)
\eqno	(6)
$$
on an interval $I \subset \R$ with locally integrable, real-valued coefficients
$m$, $l$ and $q$ (we have absorbed the spectral parameter in the latter).
If $u$ is a solution of $(6)$, then so is its (component-wise) complex
conjugate, so it is sufficient to study $\R^2$-valued solutions.
Since $(6)$ is linear, a non-trivial solution will never take the value
$\pmatrix{0 \cr 0\cr}$, so $u(x)$ traces out an absolutely continuous curve
in the punctured plane as $x$ varies.

Introducing polar coordinates in this plane by writing
$$
	u(x) = \modul{u(x)} \pmatrix{\cos \teta(x) \cr \sin \teta(x) \cr},
$$
where the absolutely continuous function $\teta$ (the {\it Pr\"ufer angle\/}
of $u$) is determined uniquely up to addition of an integer multiple of
$2 \pi$ --- or indeed $\pi$, since $-u$ is again a solution of $(6)$ ---,
a straightforward calculation yields the differential equation for $\teta$
({\it Pr\"ufer equation\/})
$$
	\teta' = m \cos 2 \teta + l \sin 2 \teta + q.
$$
This non-linear equation is equivalent to $(6)$, because one can recover
$u$ from $\teta$ by noticing that (choosing $x_0 \in I$)
$$
	\modul{u(x)} = \modul{u(x_0)}
		\exp\int_{x_0}^x (m \sin 2 \teta - l \cos 2 \teta)
\qquad	(x \in I).
$$ 

\medskip
\noindent
{\it B. Dirac systems and the factorisation method\/}

\noindent
Known to students of quantum mechanics as a trick to treat the Schr\"odinger
equation for the harmonic oscillator, the factorisation method is in fact a
tool of much wider scope for studying second-order differential equations
(see the extensive discussion in [6]).
It is based on the close link between the eigenvalue equations
$$
	A^* A u = \lambda u
\qquad	\hbox{ and }
\qquad	A A^* v = \lambda v,
$$
where $A$ is typically a first-order differential operator: a solution $u$
of the first equation gives rise to a solution $v = A u$ of the second.

	Here we remark on how this method can serve to solve eigenvalue
equations for certain one-dimensional Dirac operators, which in turn
represent a factorisation of the Sturm-Liouville operators arising as their
formal squares.

	Starting from a Dirac system $(6)$, or equivalently
$$	\eqalign{
&	u_1' = - l u_1 + (m - q) u_2,	\cr
&	u_2' = (m + q) u_1 + l u_2,	\cr }
$$
with absolutely continuous coefficients $l$, $m$, $q$, we find that $u_1$
and $u_2$ satisfy the second-order equations
$$	\eqalign{
	&u_1'' = (m^2 - q^2 + l^2 - l') u_1 + (m' - q') u_2,
\cr	&u_2'' = (m^2 - q^2 + l^2 + l') u_2 + (m' + q') u_1.
\cr	}
$$
This ordinary differential equation system decouples if $m$ and $q$ are
constant.
If we assume this in the following and set $c := m - q$, $d := m + q$ and
$a := cd$,
we are dealing with the Dirac system
$$	\eqalign{
	&u_1' = -l u_1 + c u_2
\cr	&u_2' = d u_1 + l u_2
\cr	}
$$
and corresponding second-order equations
$$	\eqalignno{
	&u_1'' = (a + l^2 - l') u_1,	& (7)
\cr	&u_2'' = (a + l^2 + l') u_2.	& (8)
\cr	}
$$
These two Sturm-Liouville equations differ only in the sign of $l'$.
Hence we can construct a chain of interlocking equation systems, along
with special solutions, in the following way.
\Thm{Theorem 3.}
Let $(a_n)_{n \in \N}$ be a sequence of real numbers with $a_1 = 0$,
and let $(l_n)_{n \in \N}$ be a sequence of real-valued absolutely continuous
functions on $I$ such that
$$
	l_{n+1}^2 + l_{n+1}' + a_{n+1} = l_n^2 - l_n' + a_n
\qquad	(n \in \N).
$$
Let $v_1(x) := \exp \int^x l_1$ $(x \in I)$, and define recursively for
$n \in \N + 1$
$$
	v_n := l_n v_{n-1} - v_{n-1}'.
$$
Then $v_n$ is a solution of $(7)$ with $a = a_n$, $l = l_n$,
and of $(8)$ with $a = a_{n+1}$, $l = l_{n+1}$.
Furthermore, for $n \in \N + 1$ and arbitrary $d_n \in \R \setminus\{0\}$,
$$
	u := \pmatrix{	-{1 \over d_n} v_n \cr
			v_{n-1} \cr	}
$$
is a solution of the Dirac system
$$	\eqalignno{
	&u_1' = -l_n u_1 + (a_n / d_n) u_2	& (9)
\cr	&u_2' = d_n u_1 + l_n u_2.
\cr	}
$$
\endThm
{\it Proof.\/}\qquad
Because of $a_1 = 0$, the first equation of the Dirac system $(9)$ with $n = 1$
decouples from the second and can easily be solved to obtain the formula for
$v_1$.
As the first component of a solution of this Dirac system, $v_1$ satisfies
equation $(7)$ with $l = l_1$, $a = a_1$. In view of the assumptions on the coefficients
$l_n$, this equation is the same as $(8)$ with $l = 2$, $a = a_2$, the
Sturm-Liouville equation
satisfied by the second component of a solution of the Dirac system for
$n = 2$. The corresponding first component of this solution can easily be
calculated from the second equation of the Dirac system. It will depend
on $d_2$, but note that this dependence can be eliminated by normalisation
as far as obtaining a solution $v_2$ of the Sturm-Liouville equation $(7)$
with $l = l_2$, $a = a_2$ is concerned.

	This process is then iterated to find the solutions for $n \in \{3,
4, \dots\}$.

\medskip
\noindent
{\bf Examples.}

\nobreak
\noindent
1.
The classical application of the factorisation method to the one-dimensional
Schr\"odinger equation for a quantum-mechanical harmonic oscillator,
$$
	-v''(x) + x^2 v(x) = \lambda v(x),
$$
is captured in the above scheme if we choose $l_n(x) := x$ $(n \in \N)$.
The hypothesis of Theorem 3 is then satisfied with $a_n := -2(n - 1)$
$(n \in \N)$, and we obtain the well-known solutions
$$
	v_n(x) = (-\Dd{x}{}+ x)^{n-1} e^{-{x^2 \over 2}}
$$
with eigenvalues $\lambda_n := 2n - 1$ $(n \in \N)$.

\noindent
2.
The case of equation (5) is treated by choosing
$l_n(x) := \kappa + (n-1) - e^{-x}$ $(n \in \N)$; we immediately find
$v_1(x) = e^{-e^{-x} + \kappa x}$ and the recursion formula of Theorem 2.

\bigskip
\noindent
{\bf References.}
\smallskipamount=2pt plus 2pt minus 1pt

\medskip
\item{1.} V. Arnold, H. Kalf and A. Schneider:
Separated Dirac operators and asymptotically constant linear systems.
{\it Math. Proc. Cambridge Phil. Soc.\/} {\bf 121} no. 1 (1997) 141--146

\smallskip
\item{2.} H. Behncke:
The Dirac equation with an anomalous magnetic moment.
{\it Math. Z.\/} {\bf 174} (1980) 213--225

\smallskip
\item{3.} H. Behncke:
The Dirac equation with an anomalous magnetic moment II. in:
{\it Proceedings of the 7th Conference on Ordinary and Partial Differential
Equations.\/} Lect. Notes in Math. {\bf 964}, Springer, Berlin 1982

\smallskip
\item{4.} H. Behncke:
Spectral properties of the Dirac equation with anomalous magnetic moment.
{\it J. Math. Phys.\/} {\bf 26} no. 10 (1985) 2556--2559

\smallskip
\item{5.} F. Gesztesy, B. Simon, B. Thaller:
On the self-adjointness of Dirac operators with anomalous magnetic moment.
{\it Proc. Amer. Math. Soc.\/} {\bf 94} (1985) 115--118

\smallskip
\item{6.} L. Infeld, T.E. Hull:
The factorisation method.
{\it Rev. Mod. Phys.\/} {\bf 23} (1951) 21--68

\smallskip
\item{7.} H. Kalf, K.M. Schmidt:
Spectral stability of the Coulomb-Dirac Hamiltonian with anomalous magnetic
moment.
{\it J. Diff. Eq.\/} (in press)

\smallskip
\item{8.} T. Kato:
{\it Perturbation theory for linear operators.\/}
Springer, Berlin 1980

\smallskip
\item{9.} K. M. Schmidt:
Relative oscillation--non-oscillation criteria for perturbed periodic
Dirac systems.
{\it J. Math. Anal. Appl.\/} {\bf 246} (2000) 591--607

\smallskip
\item{10.} K. M. Schmidt:
Eigenvalue asymptotics of perturbed periodic Dirac systems in the slow-decay
limit.
{\it Proc. Amer. Math. Soc.\/} {\bf 131} no. 4 (2002) 1205--1214

\smallskip
\item{11.} K. M. Schmidt:
Eigenvalues in gaps of perturbed periodic Dirac operators: numerical evidence.
{\it J. Comp. Appl. Math.\/} {\bf 148} (2002) 169--181

\smallskip
\item{12.} A. V. Sobolev:
Weyl asymptotics for the discrete spectrum of the perturbed Hill operator.
{\it Adv. Sov. Math.\/} {\bf 7} (1991) 159--178

\smallskip
\item{13.} J. Weidmann:
{\it Spectral theory of ordinary differential operators.\/}
Lect. Notes in Math. {\bf 1258}, Springer, Berlin 1987

\bye